\documentclass{article}
\usepackage{graphicx}
 \usepackage{mathptmx}
\usepackage{amsmath, amstext,amssymb,amsfonts}
\usepackage[english]{babel}
\usepackage{delarray}
\usepackage{mathptmx}
\usepackage{amsmath, amstext,amssymb,amsfonts}
\usepackage[english]{babel}
\usepackage{delarray}

\DeclareFontFamily{U}{mathx}{\hyphenchar\font45}
\DeclareFontShape{U}{mathx}{m}{n}{
      <5> <6> <7> <8> <9> <10>
      <10.95> <12> <14.4> <17.28> <20.74> <24.88>
      mathx10
      }{}
\DeclareSymbolFont{mathx}{U}{mathx}{m}{n}
\DeclareFontSubstitution{U}{mathx}{m}{n}
\DeclareMathAccent{\widecheck}{0}{mathx}{"71}
\DeclareMathAccent{\wideparen}{0}{mathx}{"75}

\numberwithin{equation}{section}
\newtheorem{teo}{Theorem}[section]       \newtheorem{lem}[teo]{Lemma}
     \newcommand{\R}{\mathbb R}    \newcommand{\Co}{\mathbb C}
 \newtheorem{pot}{Proof}  \newtheorem{pott}{Proof of Theorem \ref{t2} (Necessity)} \newtheorem{pottt}{Proof of Theorem \ref{t2} (Sufficiency)}

\date{}
\title{On  positive definite distributions}
\author{{\large Saulius  Norvidas} }
\date{{\footnotesize Institute of Data Science and Digital Technologies, Vilnius University, Akademijos str. 4, Vilnius LT-04812, Lithuania\\
 ({\rm{e-mail: norvidas{@}gmail.com}})}}
\begin{document}

\maketitle
 {{ {\bf Abstract}}}
We provide necessary and sufficient conditions for a tempered distribution $F\in S'(\R)$ to be positive definite. A generalized  Cauchy transform $\widetilde{F}$ of $F$ is used as a numerical continuation of $F$ to the open upper and lower complex half-planes in $\Co$. In fact, our necessary and sufficient conditions for $F$ are determined completely by  the properties of  the restriction of $\widetilde{F}$ to  the imaginary axis  in $\Co$. The main result is given in terms of completely monotonic and absolutely monotonic functions.


\section{Introduction}
\label{s1}

A function $f:\R\to\Co$ is said to be positive definite if
\begin{equation}\label{1.1}
\sum_{j, k=1}^n f(x_j-x_k)c_j{\overline{c}}_k\ge 0
\end{equation}
holds for all finite sets of complex numbers $c_1,\dots, c_n$ and points $x_1,\dots,x_n\in\R$.  A survey about positive definite functions and its generalizations can be found in \cite{9}.    The Bochner theorem states that a continuous function $f$ on $\R$ is positive definite if and only if it is the Fourier transform of a finite nonnegative measure $\mu$ on $\R$, i.e.
\[
f(x)= \hat{\mu}(x)=\int_{-\infty}^{\infty}e^{-ixt}\,d\mu(t),
\]
$x\in\R$. We define the inverse Fourier transform as
\begin{equation}\label{1.2}
\check{\mu}(\xi)=\frac1{2\pi}\int_{-\infty}^{\infty}e^{i\xi t}\,d\mu(t).
\end{equation}
In the case of integrable functions $\varphi$ on $\R$  the Fourier transform  and its inverse are defined  similarly. By our  normalization in (1.2), the  following inversion formula $\widehat{(\check{\varphi})}=\varphi$  holds  for suitable functions $\varphi$.

 For many purposes it is convenient to replace  (1.1) by its integrable analogue
\begin{equation}\label{1.3}
\int_{-\infty}^{\infty}\int_{-\infty}^{\infty} f(x-y) \varphi(x)\varphi(y)\,dx\,dy\ge 0,
\end{equation}
where $\varphi$ ranges over $L^1(\R)$ or over the space $C_c(\R)$ of continuous functions  with compact support. If $f$ is continuous, then (1.3) is equivalent to (1.1) (see [9, p. 420]).  Note that (1.3) can be rewritten in the form
\begin{equation}\label{1.4}
\int_{-\infty}^{\infty}f(x)\bigl(\varphi\ast\varphi^{\star}\bigr)(x)\,dx\ge 0,
\end{equation}
where $\varphi^{\star}(x):= {\overline{\varphi(-x)}}$, and $\ast$ denotes the convolution operation
\[
u\ast v(x)=\int_{-\infty}^{\infty}u(x-t)v(t)\,dt.
\]
The property (1.4) can be taken as the basis for defining positive definite distributions (or generalized functions). Let us recall some notion. We shall follow [10].

 Let $C^m(\R)$  denote the set of complex-valued function on $\R$ with continuous derivatives of order $m\in \mathbb{N}_0:=\mathbb{N}\cup \{0\} $, and let $C^{\infty}(\R)=\cap_{m\in \mathbb{N}_0} C^m(\R)$. The Schwartz space $S(\R)$ of test functions can be defined as the set   of  $\varphi\in C^{\infty}(\R)$ satisfying
\begin{equation}\label{1.5}
\|\varphi\|_m:=\sup_{\{x\in\R; \ 0\le p\le m\}}(1+|x|)^{m}\bigl|\varphi^{(p)}(x)\bigr|<\infty
\end{equation}
for all $m\in \mathbb{N}_0$. These  seminorms turn $S(\R)$ into a Fr\'{e}chet space. The elements of the dual space $S'(\R)$  are
called tempered distributions on $\R$.  Let $D(\R)$ denote the subspace of  $ C^{\infty}(\R)$ consisting of  functions with compact support.  The topology on  $D(\R)$ is introduced as usual (see [3] or [10]). Since the convergence in  $D(\R)$  implies convergence in $S(\R)$, it follows that  $S'(\R) \subset D'(\R)$. The elements $F$ in $ D'(\R)$ are called distributions on $\R$, and their action on test functions $\varphi\in D(\R) $ is written as $(F,\varphi)$.

  Note that a nonnegative $\sigma$-additive measure $\eta$ on  the Borel subsets of $\R$ such that $\eta(A)<\infty$ for all bounded $A$, defines   via the formula
\[
(F_{\eta}, \varphi)=\int_{-\infty}^{\infty} \varphi(x)\,d \eta(x), \ \varphi\in S(\R),
\]
the element $F_{\eta}$ in $S'(\R)$ if and only if $\eta$ is a tempered measure, i.e.
\[
\int_{-\infty}^{\infty} (1+|x|)^{-m}\,d \eta(x)<\infty
\]
for some $m\in \mathbb{N}_0$ (see [10, p. 78]).

  Let us denote by $S_m(\R)$  the completion of $S(\R)$ (as metric space) in the $m$-th norm (1.5). The following lemma gives an exact characterization of   $S_m(\R)$.

\begin{lem}\label{t1}{\rm{[10, p. 75]}}.
Let $m\in \mathbb{N}_0$. A function $u:\R\to\Co$  is an element of $S_m(\R)$ if and only if $u\in C^m(\R)$ and $x^m u^{(p)}(x)\to 0$ for $|x|\to\infty$ and all $p\le m$.
\end{lem}

 Since
\begin{equation}\label{1.6}
\|\varphi\|_0\le \|\varphi\|_1\le \|\varphi\|_2\le\dots ,\quad \varphi\in S(\R),
\end{equation}
 we see that a linear functional $F$ on $S(\R)$ is continuous  if and only if there exist  $ A>0$ and  $m\in \mathbb{N}_0$ such that
\begin{equation}\label{1.7}
|(F, \varphi)|\le A\|\varphi\|_m
\end{equation}
for all $\varphi\in S(\R)$ (see [7, p. 74]). We will call the smallest   $m\in \mathbb{N}_0$ for which (1.7) holds true for certain $A$,  the  $S'$-{\it order} of $F\in S'(\R)$. In the sequel,  $\varrho_S(F)$ denotes this order. Note   that each $F\in S'(\R)$ can be extended  to a continuous linear functional on   $S_m(\R)$ for all  $m\ge \varrho_S(F)$.  Our definition of  $\varrho_S(F)$   is different from the standard definition of order in $D'(\R)$. The usual order can be defined in various equivalent ways: for example,  $P\in D'(\R)$ is a distribution of finite order if there exists an  $m\in \mathbb{N}_0$ such that $P$ can be continued to a continuous linear functional on  $C^m_c(\R)$ [3, p. 41].  The smallest such  $m$ is called the order of $P\in D'(\R)$. We  call  this order as  the $D'$-{\it order} of $P$,  and denote it by $\varrho_D(P)$. If $F\in S'(\R)$, then it is easy to see that
\[
\varrho_D(F)\le \varrho_S(F)<\infty.
\]
 In the case where $F$ is a distribution with compact support, we have $\varrho_D(F)=\varrho_S(F)$. Suppose that $q$ is a polynomial  of degree $k$ and $F_q$ is the following regular tempered distribution
 \[
 (F_q, \varphi)=\int_{-\infty}^{\infty}q(x)\varphi(x)\,dx,\quad \varphi\in S(\R)
 \]
 associated with $q$. It is obvious that  $\varrho_D(F_q)=0$. On the other hand, by Lemma 1.1, we see that $\varrho_S(F_q)=k+2$.

 Let $a\in\R$. The operator $E_a: F\to e^{iat} F$, where $(e^{iat} F, \varphi)=(e^{iat} F_t, \varphi(t))= (F_t, e^{iat}\varphi(t))$, $\varphi\in S(\R)$,  is continuous and invertible on $S(\R)$. Moreover,  this operator acts continuously on each  $S_m(\R)$, $m\in\mathbb{N}_0$, and is not difficult to show that  $\|E_a\|_{S_m(\R)}, \|(E_a)^{-1}\|_{S_m(\R)}\le (1+|a|)^m$.  Therefore,  we get
\begin{equation}\label{1.8}
  \varrho_S(e^{iat}F)=\varrho_S(F).
\end{equation}

 Let $F\in S'(\R)$. If $m\in \mathbb{N}_0$ and  $m\ge \varrho_S(F)$, then the relation  (1.8) and  Lemma 1.1 show that for any fixed $a\in\R$,   the following generalized Cauchy transform of $F$
\begin{equation}\label{1.9}
\widetilde{F}(z)=\frac{i}{\pi}\Bigl(e^{iat}F_t, \frac{1}{(z-t)^{m+1}}\Bigr)
\end{equation}
 is well-defined for  $z\in \Co\setminus\R$. This  Cauchy transform was considered in [2, Ch. 6] (with $a=0$) for the  purposes of analytic representation of distributions $F\in {\cal{O}}'_{\alpha}$,  where  $ {\cal{O}}_{\alpha}$ is the space of  $\varphi\in C^{\infty}(\R)$ such that $\varphi^{(k)}(x)=O(|x|^{\alpha})$ as $|x|\to\infty$ for all $k\in\mathbb{N}_0$,  with a suitably  topology defined  in terms of convergent sequences.

 A distribution $F\in D'(\R)$ is said  to be positive definite if $(F, \varphi\ast\varphi^{\star})\ge 0$  for all $\varphi\in D(\R)$. The Bochner-Schwartz  theorem states that $F\in D'(\R)$   is positive definite if and only if $F$ is the Fourier transform of a nonnegative tempered measure on $\R$ [10, p. 125]. In particular, this theorem implies that any  positive definite distribution lies in $ S'(\R)$. An important problem with practical implication  is to find conditions on $F\in S'(\R)$ under which $F$  is positive definite.  There are many characterizations of continuous positive definite functions (see, for example,  [5, p.p. 70-83]). As far as we known, it is perhaps surprising that there are almost no such  results in the case of distributions.

  In this paper, we will consider the  characterization of positive definiteness for $F\in S'(\R)$ by means of some properties of the numerical function  (1.9). A similar problem for continuous  positive definite functions $f$  on $\R$ was studied in [6] by using the Poisson transform
 \[
 u_f(z)=\frac1{\pi}\int_{-\infty}^{\infty} \frac{y}{(x-t)^2+y^2}f(t)\,dt,
\]
$z=x+iy\in \Co\setminus\R$. Note that this Poisson transform is the real part of the usual Cauchy transform of $f$. Recall that a function $\omega: (a, b)\to \R$  is said to be  completely monotonic if  $\omega$ is infinitely differentiable  and $(-1)^k \omega^{(k)}(x)\ge 0$ for any $x\in(a, b)$ and all  $k\in\mathbb{N}_0$.  A function $\omega: [a, b]\to \R$  is called completely monotonic   on $[a, b]$ if $\omega$ is there continuous and completely monotonic   on $(a, b)$.
\begin{teo}\label{t1}{\rm{[6,  Theorem  2].}}
Suppose that $f: \R\to\R $ is a    bounded continuous  even function. If   $f(0)=1$, then $f$ is a characteristic function of a probability measure on $\R$ if and only if the function $ y\to u_f(0,y)$  is completely monotonic on $[0,\infty)$.
\end{teo}
  In the case where $f$ is  absolutely integrable and infinitely differentiable on $\R$, some  theorem of this kind (in other terms) has been shown by Egorov in [4].

   A function $\omega(x)$  is said to be absolutely monotonic  on $(a, b)$ if $\omega(-x)$ is completely monotonic on $(-b, -a)$.  It is obvious that such a function $\omega(x)$ can be  characterized by  the property that $\omega^{(k)}(x)\ge 0$ for  any $x\in (a,b)$ and  all  $k\in\mathbb{N}_0$.

  The main result is the following.

\begin{teo}\label{t2}
Let $F\in S'(\R)$, and suppose $n\in\mathbb{N}_0$ be  such that $2n\ge \varrho_S(F)$. Let $a_1,a_2\in\R$, where $a_1\neq a_2$, and let
\begin{equation}\label{1.10}
\widetilde{F}_j(z)=(-1)^n\frac{i}{\pi}\Bigl(e^{ia_j\,t}F, \frac{1}{(z-t)^{2n+1}}\Bigr)
\end{equation}
for $z\in \Co\setminus\R$ and $j=1,2$. Then $F$ is  positive definite if and only if:

({\bf i}) \   $ \widetilde{F}_j(iy)$, $j=1, 2$   are completely monotonic functions for $y\in (0,\infty)$;

({\bf ii}) \   $ -\widetilde{F}_j(iy)$, $j=1, 2$   are   absolutely monotonic functions for $y\in (-\infty, 0)$.
\end{teo}

It is quite possible that such a characterization of positive definite distributions on $\mathbb{R}^n$ is also valid. In this case, we could study the Cauchy-Bochner transform (instead of (1.9)) and its generalization in the tube domain $T^{\Lambda}=\mathbb{R}^n+i\Lambda\subset \mathbb{C}^n$, where  $\Lambda$ is an open cone in $\mathbb{R}^n$ (see, for example,   [10, p. 144]). On the other hand,  the class of completely monotonic functions on the $n$-dimensional algebraic structures is also studied in detail (see  [1]). We also see some technical difficulties that may arise in the case of several variables.

\section{Preliminaries and  Proofs}
 \label{s2}
Let us start  with the observation that if  $\varphi\in C^{m+1}(\R)$ and $\varphi$ has finite $(m+1)$th norm (1.5), then  by Lemma 1.1,  $\varphi$ belongs to $S_m(\R)$  but not necessarily to $S_{m+1}(\R)$. For example, we have such a case if
\[
\varphi(x)=\frac1{(i+x)^{m+1}}.
\]
\begin{lem}\label{c2}
Let $F\in S'(\R)$,  and suppose $m\in\mathbb{N}_0$ be  such that  that $m\ge \varrho_S(F)$. Let  $\varphi\in C^{\infty}(\R)$ be  a positive definite function  such that $\|\varphi\|_{m+1}<\infty$. If $F$ is positive definite, then $(F, \varphi)\ge 0$.
\end{lem}
\begin{pot}
We first claim that there exists a sequence $(\varphi_n)_{n=1}^{\infty}$ of positive definite functions $\varphi_n$ in $D(\R)$ such that $\varphi_n\to\varphi$ in the $S_m(\R)$ norm. To see this, fix any  $a_1\in D(\R)$ such that $\|a_1\|_{L^2(\R)}=1$. Let $a=a_1\ast a^{\star}_1$. Then $a\in D(\R)$, $a(0)=\|a_1\|^2_{L^2(\R)}=1$, and $\widehat{a}=|\widehat{a_1}|^2$. Therefore, $a$ is  positive definite. Define the family of continuous  positive definite functions

\[
A=\{\varphi_{\varepsilon}(x):= a(\varepsilon x)\varphi(x), \quad \varepsilon\in (0,1]\}.
\]
Since  $\|\varphi\|_{m+1}<\infty$ and $\varepsilon\le 1$, it it is easy to verify that
\[
\|\varphi_{\varepsilon}\|_{m+1}\le \sigma \|\varphi\|_{m+1}
\]
for all $\varphi_{\varepsilon}\in A$, where
\[
\sigma=\max_{0\le p\le m+1} \sum _{k=0}^p  \binom{p}{k} \|a^{(k)}\|_{L^{\infty}(\R)}<\infty.
\]
Hence $A$ is a bounded subset of $S_{m+1}(\R)$. The natural imbedding $T_{m+1}: S_{m+1}(\R)\to S_{m}(\R)$ is a compact operator [10,  p. 75]. Therefore,  $A=T_{m+1}(A)$ is totally bounded (precompact) in the $S_m(\R)$ topology. Also, there is a sequence $(\varepsilon_n)_{n=1}^{\infty}$ of positive numbers,  $\lim_{n\to\infty} \varepsilon_n=0$, and a sequence
\[
\varphi_n(x):=a(\varepsilon_n x)\varphi(x), \ n=1,2,\dots,
\]
  $\varphi_n\in A\subset  D(\R)$  such that $(\varphi_n)_{n=1}^{\infty}$  converges in the $S_m(\R)$ norm to some   $\psi\in S_m(\R)$. Since $a$  is continuous on $\R$ and $a(0)=1$, it follows that $\psi=\varphi$.

 Now the condition  $m\ge \varrho_S(F)$ guarantees  that $F$ can be extended to an element of $S'_m(\R)$. Therefore, by the definition of the Fourier transform in $S'(\R)$, we get
\[
(F, \varphi)=\lim_{n\to\infty} (F, \varphi_n)=\lim_{n\to\infty} (\hat{F}, \check{\varphi}_n).
\]
 According to the  Bochner and the Bochner--Schwartz theorems, we have that $\check{\varphi}_n$ and $\hat{F}$ are nonnegative elements of $S(\R)$ and $S'(\R)$,  respectively. This means, in particular, that   $\hat{F}$  is a nonnegative tempered measure $\mu$ on $\R$. Hence
 \[
 (\hat{F}, \check{\varphi}_n)=\int_{-\infty}^{\infty} \check{\varphi}_n(x)\, d\mu(x)\ge 0
 \]
  and thus $(F, \varphi)\ge 0$.

\end{pot}

 We need the following technically lemma; for completeness, we also give its proof.
\begin{lem}\label{c3}
Let $F\in S'(\R)$,  $a\in\R$, and let $m\ge \varrho_S(F)$. The function
\[
\widetilde{F}(z)=\Bigl(e^{iat}F_t, \frac{1}{(z-t)^{m+1}}\Bigr)
\]
 is analytic in $\Co\setminus\R$. Moreover, for any  $s\in\mathbb{N}$,
\begin{equation}\label{2.1}
\frac{d^s}{d z^s}\,\widetilde{F}(z)=\biggl(e^{iat} F_t, \frac{d^s}{dz^s}\biggl(\frac1{(z-t)^{m+1}}\biggr)\biggr).
\end{equation}
\end{lem}
\begin{pot}
Fix $z\in\Co\setminus\R$. Take any $h\in \Co\setminus\{0\}$  such that $z-h\in\Co\setminus\R$. Since
\begin{gather}
\frac1{h}\Bigl[\widetilde{F}(z+h)-\widetilde{F}(z)\Bigr]=\biggl(e^{iat}F, \frac{1}{h}\biggl[\frac1{(z+h-t)^{m+1}}-\frac1{(z-t)^{m+1}}\biggr]\biggr)
\nonumber
\end{gather}
and $\varrho_S(e^{iat} F)= \varrho_S( F)\le m$,  we see that $\widetilde{F}$ is analytic at $z$ if and only if
\begin{equation}\label{2.2}
\frac{1}{h}\biggl[\frac1{(z+h-t)^{m+1}}-\frac1{(z-t)^{m+1}}\biggr]\to \frac{-(m+1)}{(z-t)^{m+2}}
\end{equation}
as $h\to 0$ in the sense of convergence in $S_m(\R)$, i.e. as functions of $t\in \R$ in the $m$-th norm (1.5). Define
\[
E_m(t)=\frac{1}{h}\biggl[\frac1{(z+h-t)^{m+1}}-\frac1{(z-t)^{m+1}}\biggr], \ t\in\R.
\]
Then $E_m\in S_{m+1}(\R)\subset S_m(\R)$ and
\[
\frac{d^p}{dt^p} E_m(t)=\frac{(m+p)!}{m!} E_{m+p}(t)
\]
for each $p\in\mathbb{N}$. Hence
\[
\biggl\|E_m(t)+\frac{m+1}{(z-t)^{m+2}}\biggr\|_{S_m(\R)}=\frac{(m+p)!}{m!}\sup_{t\in\R;\ 0\le p\le m}\biggl|(1+|t|)^m\biggl[E_{m+p}(t)+\frac{m+p+1}{(z-t)^{m+p+2}}\biggr]\biggr|.
\]
Now a straightforward calculation leads to the identity
\begin{gather}\label{2.3}
\biggl\|E_m(t)+\frac{m+1}{(z-t)^{m+2}}\biggr\|_{S_m(\R)}
\nonumber      \\
=\frac{(m+p)!}{m!}\sup_{t\in\R;\ 0\le p\le m}\biggl(\frac{(1+|t|)^m}{|z+h-t|^{m+p+1}|z-t|^2}\biggl|\sum_{j=1}^{m+p+1} a_j(m,p)\frac{h^j}{(z-t)^{j-1}}\biggr|\biggr),
\end{gather}
where the coefficients $a_j(m,p)$ depend only on $m$ and $p$. Clearly, if $h\to 0$, then using (2.3) we obtain (2.2). This shows  (2.1)  in the  case $s=1$. Therefore,  $\widetilde{F}$ is analytic in $\Co\setminus\R$. The cases of higher $s\in \mathbb{N}$  can be obtained  by induction. This is a standard procedure   and we omit the details.
\end{pot}
 For the sake of clarity, we divide the proof of the main result in two parts.

 \begin{pott}
  Fix any $j=1,2$. Using  Lemma 2.2 and having in mind that $(-1)^ni=(i)^{2n+1}$, we get
\begin{gather}\label{2.4}
\frac{d^s}{dy^s}\widetilde{F}_j(iy)=\frac1{\pi} \biggl(e^{ia_jt} F, \frac{\partial^s}{\partial y^s}\biggl(\frac1{(y+it)^{2n+1}}\biggr)\biggr)=
\frac{(-1)^s}{\pi}\frac{(2n+s)!}{(2n)!} \biggl(e^{ia_jt} F, \frac{1}{(y+it)^{2n+s+1}}\biggr)\nonumber      \\
=
\frac{(-1)^s}{\pi}\frac{(2n+s)!}{(2n)!} \biggl( F, \frac{e^{ia_jt}}{(y+it)^{2n+s+1}}\biggr)
\end{gather}
for all $s\in\mathbb{N}_0$.

 Suppose that $y>0$. By direct calculation, we have
\[
\int_0^{\infty} x^{2n+s} e^{-yx} e^{-ixt}\,dx= \frac{(2n+s)!}{(y+it)^{2n+s+1}}.
\]
 Bochner's theorem shows that the right side of this equality is  positive definite as a function of $t$ in $\R$. Furthermore, since the product of positive definite functions is also positive definite, it follows that the function
\[
\varphi_j(t)=\frac{e^{ia_j t}}{(y+it)^{2n+s+1}}
\]
 satisfies the conditions of Lemma 2.1 for $m=2n+s$. Therefore, using (2.4) we obtain that
\[
(-1)^s\frac{d^s}{dy^s}\widetilde{F}_j(iy)=\frac{(2n+s)!}{(2n)!\pi}(F, \varphi_j)\ge0
\]
for all $s\in\mathbb{N}_0$. This shows that the function $y\to \widetilde{F}_j(iy)$ are completely monotonic on $(0,\infty)$.

 Let now $y<0$. Since
\[
\int_{-\infty}^0 x^{2n+s} e^{-yx} e^{-ixt}\,dx= (-1)^{s+1}\frac{(2n+s)!}{(y+it)^{2n+s+1}},
\]
we see that  in this case
\[
\varphi_j(t)=(-1)^s\frac{e^{ia_j t}}{(y+it)^{2n+s+1}}
\]
 satisfies the conditions of Lemma 2.1 for $m=2n+s$. Again, by  (2.4), we get
\[
-\frac{d^s}{dy^s}\widetilde{F}_j(iy)=\frac{(2n+s)!}{(2n)!\pi}(F, \varphi_j)\ge0
\]
for all $s\in\mathbb{N}_0$. This means that $y\to -\widetilde{F}_j(iy)$ is absolute monotonic on $(0,\infty)$.  The necessity of Theorem 1.3 is proved.
 \end{pott}
 To prove the next lemma we require the concept of the Bochner integral.  For precise details we refer to [8, Ch. 1].

 Let  $(X,{\cal{B}}, m)$ be a measurable space and let $E$ be a Banach space. A function $\zeta: X\to E$ is called finite-valued if
\[
\zeta(x)=\sum_{k=1}^{n}a_k\chi_{E_{k}}(x),
\]
where $a_k\in E$, $k=1,\dots,n$,   $E_k$ are disjoint subsets of $X$ with $m(E_k)<\infty$,  and  $\chi_{E_k}$ is the indicator function of $E_k$. A function $\zeta(x): X\to E$ is said to be Bochner $m$-integrable, if there exists a sequence of finite-valued functions $\bigl(\zeta_n(x)\bigr)_{n=1}^{\infty}$ such that $\lim_{n\to\infty} \|\zeta_n(x)- \zeta(x)\|_{E}=0$ almost everywhere with respect to $m$ on $X$ and
\[
\lim_{n\to\infty}\int_X\bigl\|\zeta_n(x)-\zeta(x)\bigr\|_E\,dm(x)=0.
\]
For any $I\in {\cal{B}}$, the Bochner integral of $\zeta(x)$ on $I$ is defined by
\[
\int_I \zeta(x)\,dm(x)=\lim_{n\to\infty}\int_X\zeta(x)\chi_I(x)\,dm(x).
\]
The Bochner theorem for such vector-valued integrals states that a measurable function $\zeta(x)$ is Bochner integrable if and only if $\|\zeta(x)\|_E$ is $m$-integrable.  A corollary of this theorem says that
\begin{equation}\label{2.5}
F\Bigl(\int_I \zeta(x)\,dm(x)\Bigr)=\int_I F(\zeta(x))\,dm(x)
\end{equation}
 for any $F\in E'$. Recall that a function  $ \zeta: X\to E$ is called  weakly measurable if,   for any $F\in E'$,  the numerical function $F(\zeta(x))$ of $x\in X$ is $\cal{B}$-measurable.  Finally we note that if $E$ is a separable Banach space, then  $ \zeta$ is measurable if and only if $ \zeta$ is weakly measurable.

\begin{lem}\label{c4}
Let $u\in C(\R)\cap L^1(\R)$, $m\in \mathbb{N}_0$, and let $y\in \R\setminus\{0\}$. Then
\begin{equation}\label{2.6}
\widetilde{u}(t,y)=\frac{i}{\pi}\int_{-\infty}^{\infty}\frac1{(x+iy-t)^{m+1}} u(x)\,dx
\end{equation}
  is as function of $t\in\R$   an element of $S_m(\R)$. Furthermore, if $F\in S'_m(\R)$, then
\end{lem}
\begin{equation}\label{2.7}
(F,\widetilde{u})= \frac{i}{\pi}\int_{-\infty}^{\infty}\biggl(F_t, \frac1{(x+iy-t)^{m+1}}\biggr) u(x)\,dx.
\end{equation}

 \begin{pot}
 Since $u\in L^1(\R)$, we have that  $u$ defines by
\[
(U, \varphi)=\int_{-\infty}^{\infty}  u(x)\varphi(x)\,dx,\quad \varphi\in S_m(\R),
\]
a functional $U\in S'_m(\R)$. Then
\[
\widetilde{u}(t,y)=\frac{i}{\pi}\biggl( U_x, \frac1{(x+iy-t)^{m+1}} \biggr),
\]
$t\in\R$. In a manner similar to the proof of Lemma 2.2, we obtain that $\widetilde{u}$ is infinitely differentiable in $t$ and
\[
\frac{d^s}{dt^s}\widetilde{u}(t,y)=\frac{i}{\pi}\biggl(U_x, \frac{\partial^s}{\partial t^s}\biggl(\frac1{(x+iy-t)^{m+1}} \biggr)\biggr)
\]
for any $s\in \mathbb{N}$. Now it is easy to verify that $\|\widetilde{u}\|_{m+1}<\infty$. This shows that $\widetilde{u}\in S_m(\R)$.

 We claim that the integral (2.6) can be considered as a Bochner integral on $\R$
\begin{equation}\label{2.8}
\int_{-\infty}^{\infty}\zeta(t, x)\,dx,
\end{equation}
with the  usual Lebesgue measure  on $\R$  and functions
\begin{equation}\label{2.9}
 \zeta(t, x)=\frac{i}{\pi}\frac{u(x)}{(x+iy-t)^{m+1}}
\end{equation}
 with  values in $S_m(\R)$. Indeed,  if $F\in S'_m(\R)$, then
 \[
 F\Bigl(\zeta(x,t)\Bigr)=\frac{i}{\pi}\biggl(F_t, \frac1{(x+iy-t)^{m+1}}\biggr)\,u(x).
 \]
The condition $u\in C(\R)$ and  Lemma 2.2 yield that this function is continuous as a function of $x\in\R$. Hence $\zeta$ is weakly measurable. But, since $S_m(\R)$ is the completion of the  separable space $S(\R)$ in $m$th norm (1.5), we get that $S_m(\R)$ is also separable. Therefore, the function $t\to  \zeta(x,t)$, which acts from $\R$ into $ S_m(\R)$,  is measurable. Further, it is obvious that the  norm
\[
\sup_{\{t\in\R; \ 0\le p\le m\}}(1+|t|)^{m}\biggl|\frac{\partial^m}{\partial t^m}\biggl(\frac1{(x+iy-t)^{m+1}}\biggr)\biggr|<\infty
\]
of the function
\[
t\to \frac1{(x+iy-t)^{m+1}}
\]
 depends only on $y$ and not on $x$. Denote this norm by ${\theta}_{y}$. Then
\[
\int_{-\infty}^{\infty}\|\zeta(t, x)\|_{S_m(\R)}\,dx=\frac{\theta_y}{\pi}\|u\|_{L^1(\R)}<\infty.
\]
By the Bochner theorem for  vector-valued integrals, we obtain  that (2.6) exists as  as a Bochner integral on $\R$.

  Next, it is obvious that any  linear functional $F_a(\varphi)=\varphi(a)$, where $a\in\R$,  is continuous on $S_m(\R)$. Hence,  (2.6) gives the same element $\tilde{u}$ in $S_m(\R)$ as the Bochner integral and as  the usual numerical Lebesgue integral  depending on the parameters  $t$ and $y$. Finally, the property (2.5) shows (2.7) and the lemma is proved.
 \end{pot}

 In the following lemma we give a Plemelj type relation in  $S_m(\R)$. The Plemelj formulas for the other spaces of test functions and distributions are discussed in [2, Ch. 5-6].

\begin{lem}\label{c4}
Let $u\in D(\R)$, $m\in \mathbb{N}_0$, $y>0$, and let $\widetilde{u}(t,y)$ be the function defined by (2.6). If $y\to 0$, then
\begin{equation}\label{2.10}
\frac12\Bigl( \widetilde{u}(t,y)- \widetilde{u}(t,-y)\Bigr)\to \frac{1}{m!} u^{(m)}(t)
\end{equation}
 in the $S_m(\R)$ norm.
\end{lem}
 \begin{pot}
Fix $\varepsilon>0$. First, we claim that there is an $A>0$ such that
\begin{equation}\label{2.11}
I_1:=\sup_{\{|t|\ge A; 0\le p\le m\}}(1+|t|)^m\Bigl| \frac{d^p}{dt^p}\biggl(\frac12\Bigl(\widetilde{u}(t,y)- \widetilde{u}(t,-y)\Bigr)- \frac{1}{m!} u^{(m)}(t)\biggr)\Bigr|<\varepsilon
\end{equation}
for all $y>0$. Indeed, if $u$ is supported on $[-a, a]$, $a>0$, then for every $A>a$ we have
\begin{gather}
I_1=\sup_{\{|t|\ge A; 0\le p\le m\}}\frac{(1+|t|)^m}{2}\Bigl|  \frac{d^p}{dt^p}\Bigl(\widetilde{u}(t,y)- \widetilde{u}(t,-y)\Bigr)\Bigr|\nonumber\\
\le
\frac{(m+p)!}{\pi m!}\sup_{\{|t|\ge A; 0\le p\le m\}}\int_{-a}^a\frac{(1+|t|)^m}{|x+iy-t|^{m+p+1}}|u(x)|\,dx\nonumber\\
\le
\frac{(m+p)!}{\pi m!}\sup_{\{|t|\ge A; 0\le p\le m\}}\int_{-a}^a\frac{(1+|t|)^m}{|x-t|^{m+p+1}}|u(x)|\,dx
\nonumber
\end{gather}
Now it is easy to see that there exists a number $A_0>0$ such that (2.11) is satisfied for each $A>A_0$ and all $y>0$. Fix any such $A$ to the end of the proof.

 Second, we will show that there is an $y_0>0$ such that if $y\in(0,y_0)$, then
\begin{equation}\label{2.12}
\sup_{|t|\le A; 0\le p\le m}(1+|t|)^m\Bigl| \frac{d^p}{dt^p}\biggl(\frac12\Bigl(\widetilde{u}(t,y)- \widetilde{u}(t,-y)\Bigr)- \frac{1}{m!} u^{(m)}(t)\biggr)\Bigr|<\varepsilon
\end{equation}
To this end, we will transform  $ \widetilde{u}(t,y)- \widetilde{u}(t,-y)$. Denote
\[
\alpha(x)=\frac1{x-iy-t}-\frac1{x+iy-t}={\overline{\biggl(\frac1{x+iy-t}-\frac1{x-iy-t}\biggr)}}.
\]
Then $\alpha\in L^1(\R)$ and
\[
\hat{\alpha}(\xi)=\int_{-\infty}^{\infty}\alpha(x)e^{-i x\xi}\,dx=2\pi i e^{-y|\xi|}e^{-i\xi t}.
\]
Integrating by parts in (2.6) and applying the Parseval formula
\[
\int_{-\infty}^{\infty}{\overline{\alpha(x)}}\beta(x)\,dt=\frac1{2\pi}\int_{-\infty}^{\infty} {\overline{\widehat{\alpha}(\xi)}}\widehat{\beta}(\xi)\,d\xi
\]
for $\beta(t)=u^{(m)}(x)$, we have
\begin{gather}
\widetilde{u}(t,y)- \widetilde{u}(t,-y)=\frac{i}{\pi m!}\int_{-\infty}^{\infty}\biggl(\frac1{x+iy-t}-\frac1{x-iy-t}\biggr)u^{(m)}(x)\,dx\\ \nonumber
=\frac{1}{\pi m!}\int_{-\infty}^{\infty}e^{-y|\xi|}e^{i\xi t}\widehat{u^{(m)}}(\xi)\,d\xi.
\end{gather}
Since $u\in D(\R)$, we may conclude that
\[
u^{(m)}(t)=\frac1{2\pi}\int_{-\infty}^{\infty}e^{i\xi t}\widehat{u^{(m)}}(\xi)\,d\xi.
\]
Now, by (2.13),
\[
\frac12\Bigl(\widetilde{u}(t,y)- \widetilde{u}(t,-y)\Bigr)- \frac{1}{m!} u^{(m)}(t)=\frac{1}{2\pi m!}\int_{-\infty}^{\infty}\Bigl(e^{-y|\xi|}-1\Bigr)e^{i\xi t}\widehat{u^{(m)}}(\xi)\,d\xi.
\]
Hence,
\begin{gather}
\frac{d^p}{dt^p}\biggl(\frac12\Bigl(\widetilde{u}(t,y)- \widetilde{u}(t,-y)\Bigr)- \frac{1}{m!} u^{(m)}(t)\biggr)\\ \nonumber
=\frac{1}{2\pi m!}\int_{-\infty}^{\infty}\Bigl(e^{-y|\xi|}-1\Bigr)(i\xi)^pe^{i\xi t}\widehat{u^{(m)}}(\xi)\,d\xi.
\end{gather}
Here the term $\chi(\xi)=(i\xi)^pe^{i\xi t}\widehat{u^{(m)}}(\xi)$ belongs to $S(\R)$. Therefore, with  the relation (2.14)  in mind, we can deduce  that there exists $y_0>0$ such that the inequality (2.12) is satisfied for all $y\in(0,y_0)$.

 Finally, combining (2.11) with (2.12), we obtain that
\[
\Bigl\|\frac12\Bigl( \widetilde{u}(t,y)- \widetilde{u}(t,-y)\Bigr)\to \frac{1}{m!} u^{(m)}(t)\Bigr\|_{S_m(\R)}<\varepsilon
\]
for all $y\in(0,y_0)$.  This proves (2.10).
 \end{pot}
 \begin{pottt}
 Fix any $j=1,2$,  and suppose that $\widetilde{F}_j(iy)$ is completely monotonic for $y\in(0,\infty)$. This means that for any $\tau>0$, the function $u_{\tau}(iy):=\widetilde{F}_j(i(\tau+y))$ is continuous and completely monotonic on $[0,\infty)$. By the Bernstein-Widder theorem (see [1, p. 135]), there is a nonnegative finite measure $\mu_{\tau}$ on $[0,\infty)$ such that
\[
 u_{\tau}(iy)=\int_0^{\infty}e^{-y\,t}\,d\mu_{\tau}(t)
 \]
 for $y\in [0,\infty) $. This function  can be continued analytically to the open upper complex half-plane $\Co^{+}$ as the Laplace transform of a finite measure by
\begin{equation}\label{2.15}
 u_{\tau}(z)=\int_0^{\infty}e^{iz\,t}\,d\mu_{\tau}(t), \quad z\in \Co^{+}.
\end{equation}
Since  $\mu_{\tau}$ is  nonnegative,    Bochner's theorem show that for any fixed $y\ge 0$ in (2.15), the function $x\to u_{\tau}(x+iy)$ is continuous positive definite for $x\in\R$. On the other hand, according to Lemma 2.2, the function $\widetilde{F}_j$ is analytic on $\Co\setminus\R$. The uniqueness theorem for the analytic functions (2.15) and $\widetilde{F}_j$ yields $\widetilde{F}_j(z)=u_{\tau}(z-i\tau)$ for all $\Im z\ge \tau$. Since $\tau$ is arbitrary positive, we see, in particular,  that for any fixed $y>0$, the function $x\to \widetilde{F}_j(x+iy)$ is also continuous positive definite for $x\in\R$. This means that
\begin{equation}\label{2.16}
\int_{-\infty}^{\infty}\widetilde{F}_j(x+iy) u(x)\,d x \ge 0
\end{equation}
holds for any $y>0$ and all positive definite functions $u\in D(\R)$.

 Suppose now that $-\widetilde{F}_j(iy)$ is absolutely monotonic for $y\in(-\infty, 0)$. Then the function $-\widetilde{F}_j(-iy)$ is completely monotonic on $(0, \infty)$. Using the same argument as before, we have
\begin{equation}\label{2.17}
\int_{-\infty}^{\infty}\widetilde{F}_j(x-iy) u(x)\,d x \le 0
\end{equation}
for any $y>0$ and all positive definite $u\in D(\R)$.

  Let $y>0$ and suppose that $u\in D(\R)$. Having in mind the definition (1.10) of $\widetilde{F}_j$ and Lemma 2.3 in the case $m=2n$, we get
\[
\int_{-\infty}^{\infty}\Bigl(\widetilde{F}_j(x+iy) -\widetilde{F}_j(x-iy)\Bigr) u(x)\,d x =(-1)^n\Bigl(e^{ia_j\,t}F, \tilde{u}(t, y)-\tilde{u}(t, -y)\Bigr).
\]
If a function  $u$ is, in addition, positive definite on $\R$, then  combining the last relation with  (2.10), (2.16), and (2.17), we obtain
\begin{equation}\label{2.18}
(-1)^n\frac{2}{(2n)!}\Bigl(e^{ia_j\,t}F, u^{(2n)}(t)\Bigr)=\lim_{y\to 0}\int_{-\infty}^{\infty}\Bigl(\widetilde{F}_j(x+iy) -\widetilde{F}_j(x-iy)\Bigr) u(x)\,d x\ge 0.
\end{equation}
Note that by the definition of  derivatives in $S'(\R)$,
 \[
 \Bigl(e^{ia_j\,t}F, u^{(2n)}(t)\Bigr)=\Bigl(\Bigl(e^{ia_j\,t}F\Bigr)^{(2n)}, u(t)\Bigr).
\]
If we combining this with  (2.18), we conclude that  $(-1)^n\Bigl(e^{ia_j\,t}F\Bigr)^{(2n)}$ is a positive definite  distribution in $S'(\R)$. Therefore, its Fourier transform
\[
\Bigl((-1)^n\Bigl(e^{ia_j\,t}F\Bigr)^{(2n)}\Bigr)^{\widehat{}}= \xi^{2n}\widehat{F}(\xi-a_j)
\]
is a nonnegative tempered measure on $\R$, by the Bochner-Schwartz theorem. In particular, we see that $\widehat{F}$ is a nonnegative distribution on $\R\setminus\{-a_j\}$.

Now, if we recall that $a_1\neq a_2$,  we  conclude   that $\widehat{F}$ is a nonnegative distribution on the whole $\R$. Therefore, $\widehat{F}$ is a nonnegative tempered measure on $\R$ (see [10, p. 17]). Finally, the Bochner-Schwartz theorem shows that $F$ is positive definite. The theorem is completely proved.
\end{pottt}
\vspace{5mm}
{\bf{Acknowledgement}}

 The author thanks the referee for pointing out several mistakes  and  making a few other remarks which  improved  the exposition.

This research was funded by a grant (No. MIP--47/2010) from the Research Council of Lithuania.
\vspace{5mm}

\bibliographystyle{elsarticle-harv}

\begin{thebibliography}{00}

\bibitem[1]{1} C. Berg, J.P.R. Christensen, P. Ressel,   Harmonic Analysis on Semigroups, Springer-Verlag, New York, 1984.
\bibitem[2]{2} H.J. Bremermann,  Distributions, Complex Variables, and Fourier Transforms, Addison-Wesley, Reading, Mass. etc., 1965.
\bibitem[3]{3} J.J. Duistermaat, J.A.C. Kolk,  Distributions: Theory and Applications, Birkh\"{a}user, New York, 2010.
\bibitem[4]{4} A.V. Egorov, On the theory of characteristic functions (English, Russian original), Russ. Math. Surv.  59(3) 567-568 (2004); translation from Usp. Mat. Nauk  59(3) (2004)  167-168.
\bibitem[5]{5} E. Lukacs,   Characteristic Functions.  2nd ed.,   Hafner Publishing Co., New York, 1970.
\bibitem[6]{6} S. Norvidas,  On harmonic continuation of characteristic functions,  Lith. Math. J.  50(4) (2010)  418--425.
\bibitem[7]{7} H.H. Schaefer,  Topological Vector Spaces. 2nd ed.,  Graduate Texts in Mathematics. 3, Springer-Verlag, New York etc., 1999.
\bibitem[8]{8} \u{S}. Schwabik, Y. Guoju,  Topics in Banach Space Integration, World Scientific, Singapore etc., 2005.
\bibitem[9]{9} J. Stewart, Positive definite functions and its generalizations, an historical  survey, Rocky Mountain J. Math.  6 (1976) 409-434.
\bibitem[10]{10}  V.S. Vladimirov, Methods of the Theory of Generalized Functions, Taylor $\&$ Francis, London, 2002.

 \end{thebibliography}

\end{document}